\documentclass[leqno,12pt]{article}
\usepackage{amssymb,amsfonts}
\usepackage{amsmath,latexsym}
\usepackage{french}
\usepackage{amstext,epic,eepic,epsf,pslatex}
\usepackage{graphicx}

\newcommand{\bepr}{{\em Proof} } 
\newcommand{\enpr}{\hfill \rule{.5em}{.5em}}

\newcommand{\R}{{\mathbb R}}

\newcommand{\A}{{\mathbb A}}

\newcommand{\Tr}{\hbox{Tr\,}}

\def\Xint#1{\mathchoice
{\XXint\displaystyle\textstyle{#1}}%
{\XXint\textstyle\scriptstyle{#1}}%
{\XXint\scriptstyle\scriptscriptstyle{#1}}%
{\XXint\scriptscriptstyle\scriptscriptstyle{#1}}%
\!\int}
\def\XXint#1#2#3{{\setbox0=\hbox{$#1{#2#3}{\int}$ }
\vcenter{\hbox{$#2#3$ }}\kern-.6\wd0}}

\newtheorem{defin}{Definition}[section] 
\newtheorem{prop}{Proposition}[section] 
\newtheorem{thm}{Theorem}[section]

\newtheorem{rque}{Remark}[section]

\begin{document}

\title{Mixed determinants, Compensated Integrability and new {\em a priori} estimates in Gas dynamics}

\author{Denis Serre \\ \'Ecole Normale Sup\'erieure de Lyon\thanks{U.M.P.A., UMR CNRS--ENSL \# 5669. 46 all\'ee d'Italie, 69364 Lyon cedex 07. France. {\tt denis.serre@ens-lyon.fr}}}

\date{}

\maketitle

\centerline{\em To Constantine Dafermos, with admiration, friendship and gratitude.}

\bigskip

\begin{abstract} 
We extend the scope of our recent Compensated Integrability theory, by exploiting the multi-linearity of the determinant map over ${\bf Sym}_n(\R)$. This allows us to establish new {\em a priori} estimates for inviscid gases flowing in the whole space $\R^d$. Notably, we estimate the defect measure (Boltzman equation) or weighted spacial correlations of the velocity field (Euler system). As usual, our bounds involve only the total mass and energy of the flow.
\end{abstract}

\paragraph{Keywords:} Div-BV tensors, mixed determinant, Schur complement, gas dynamics.

\paragraph{2020 Mathematics subject classification:} 35Q31, 46N60, 46T99, 76N15

\paragraph{Notations:} The unit sphere of $\R^n$ is $S_{n-1}$, whose area is $|S_{n-1}|$. For $1\le j\le n$ and $x\in\R^n$, $\hat x_j\in\R^{n-1}$ denotes the vector in which the $j$th coordinate is omitted.

The exponent $T$ indicates transposition. Given a vector $Z$, we set $Z\otimes Z:=ZZ^T$ the symmetric, rank-one matrix of entries $z_iz_j$. The cofactor matrix of an $n\times n$ matrix $M$ is $\widehat M$. We recall the identities $M^T\widehat{M}=\widehat{M}M^T=(\det M)I_n$ and $\det\widehat{M}=(\det M)^{n-1}$. The space of $n\times n$ symmetric matrices with real entries being ${\bf Sym}_n$, the cone of positive semi-definite ones is ${\bf Sym}_n^+$. Its interior is ${\bf SPD}_n$. The inequality $A\succ B$ between symmetric matrices, means that $A-B\in{\bf Sym}_n^+$.

An inequality $F(X)\le_nG(X)$ means that there exists a constant $C(n)$, depending only upon the ambiant dimension $n$, such that every $X$ under consideration satisfies $F(X)\le C(n)\cdot G(X)$. The total mass of a finite Radon measure $\mu$ is $\|\mu\|_{\cal M}$. If $\vec\mu$ is a vector of Radon measures, its Euclidian norm $|\vec\mu|$ is still a Radon measure~; we again write $\|\vec\mu\|_{\cal M}$ for the mass of $|\vec\mu|$.


\section{Introduction}

We recall that a {\em symmetric tensor} over $\R^n$ is a symmetric $n\times n$ matrix $A$, whose entries are distributions, $a_{ij}\in{\cal D}(\R^n)$. The Divergence (with a capital letter in this context) operator associates with $A$ a vector of distributions,
$$({\rm Div}\,A)_i=\sum_{j=1}^n\partial_ja_{ij},\qquad\forall i\in[\![1,n]\!].$$
Recall that if $A$ is positive semi-definite, then the $a_{ij}$'s are Radon measures. Since they are absolutely  continuous with respect to $\Tr\,A$, and because $\det^{\frac1n}:{\bf Sym}_n^+\to\R$ is positively homogeneous of order $1$, we can define unambiguously the Radon measure $(\det A)^{\frac1n}$.

\begin{defin}
The tensor $A$ is said {\em Div-BV} over $\R^n$ if its entries $a_{ij}$, as well as the coordinates $({\rm Div}\,A)_i$, are Radon measures with finite total masses. 
\end{defin}

The fundamental statement of Compensated Integrability, established in \cite{Ser_DPT,Ser_JMPA} reads as follows. 
\begin{thm}[\cite{Ser_DPT,Ser_JMPA}]\label{th:fund}
Let $A$ be a positive semi-definite Div-BV tensor over $\R^n$. Then $(\det A)^{\frac1n}\in L^{\frac n{n-1}}(\R^n)$ and
$$\int_{\R^n}(\det A)^{\frac1{n-1}}dx\le c_n\|{\rm Div}\,A\|_{\cal M}^{\frac n{n-1}}.$$
\end{thm}
This inequality is sharp, reducing to the isoperimetric inequality when $A=\chi_\Omega I_n$.

\bigskip

Theorem \ref{th:fund} implies several variants that have been described in some other papers of ours. Of particular interest is a version taylored for Cauchy problems in the whole space $\R^d$. There, a tensor is given in the strip $Q_T=(0,T)\times\R^d$, which is often  Div-free. We set $n=1+d$ and $x=(t,y)$ where $y$ is the space variable. The tensor is extended by $0_n$ away from $Q_T$. This extension is Div-BV provided that the traces of the first column (that is the normal traces) at $t=0^+$ and $t=T^-$ are finite measures. 

For the most important example of inviscid gas dynamics, the tensor is given in $Q_T$ by
$$\begin{pmatrix} \rho & \rho u^T \\ \rho u & \rho u\otimes u +pI_d \end{pmatrix}$$
where $\rho\ge0$ is the mass density, $p\ge0$ the pressure and $u$ is the fluid velocity. The specific internal energy $e\ge0$ is given in terms of $\rho,p$ through an equation of state, typically $p=(\gamma-1)\rho e$ for a perfect gas, $\gamma>1$ being the adiabatic constant. The trace condition is ensured by the assumption of finite total mass and total mechanical energy. The application of C.I. has thus led to the following estimate of the internal variables:
\begin{thm}[\cite{Ser_DPT,Ser_JMPA}]\label{th:pGD}
Consider an admissible inviscid gas flow in $Q_T$ of finite total mass and energy,
$$M:=\int_{\R^d}\rho(0,y)\,dy<+\infty,\qquad E_0:=\int_{\R^d}\left(\frac12\,\rho|u|^2+\rho e\right)(0,y)\,dy<+\infty.$$
Then 
\begin{equation}\label{eq:pGD}
\int_{Q_T}\rho^{\frac1d}p(t,y)\,dy\,dt\le_dM^{\frac1d}\sqrt{ME_0\,}\,.
\end{equation}
\end{thm}
We recall that since the right-hand side is independent from $T$, (\ref{eq:pGD}) is valid when $T=+\infty$. The same remark holds true for all the statements below, unless otherwise specified.

The admissibility assumed in Theorem \ref{th:pGD} is the fact the total mass is conserved and the total energy at positive times is bounded by that at time $t=0$~:
\begin{equation}\label{eq:totME}
\forall t\in(0,T)\qquad\int_{\R^d}\rho(t,y)\,dy=M,\quad\int_{\R^d}\left(\frac12\,\rho|u|^2+\rho e\right)(t,y)\,dy\le E_0.
\end{equation}
The strength of Estimate (\ref{eq:pGD}), when compared with (\ref{eq:totME}.b),
is the presence of an extra factor $\rho^{\frac1d}$ in the integrand. The price to pay is the replacement of a supremum over $t\in(0,T)$ by a time integral. Somehow, (\ref{eq:pGD}) plays the role of a Strichartz inequality for Gas dynamics.

The flaw of Theorem \ref{th:pGD} is of course the lack of a corresponding estimate for the velocity field. This is unsatisfactory since the energy estimate (\ref{eq:totME}.b) suggests that $\rho|u|^2$ and $\rho e$ (proportional to $p$) are on the par, having the same physical dimension. We should expect a space-time estimate of $\rho^{1+\frac1d}|u|^2$, as a counterpart of (\ref{eq:pGD}).

\bigskip

The purpose of this paper is thus to establish complementary estimates for Gas dynamics. They have several features in common with (\ref{eq:pGD}). First, they concern space-time integrals. Second, the integrands contain an extra factor not present in the energy density. Third, the bounds are still given in terms of $M$ and $E_0$. At last these estimates are invariant under Galilean transformations. Somehow, all these estimate might be viewed as {\em fundamental}, or universal~; in particular, they do not depend upon the particular equation of state. 

Most importantly, we fill the gap mentionned above, by handling the velocity field. Mind however that since we are only able to treat Galilean invariant quantities, the velocity enters in the integrals through correlations
$${\rm Cor}(u_0,\ldots,u_d)=\det\begin{pmatrix} 1 & \cdots & 1 \\ u_0 & \cdots & u_d \end{pmatrix},$$
where $u_j(t,y):=u(t,y+h_j)$ and the constant shifts $h_j\in\R^d$ are parameters. Notice that the quantity above is $d!$ times the (signed) $d$-volume of the simplex spanned by $(u_0,\ldots,u_d)$. Thus the estimate below (where of course $\rho_j(t,y)=\rho(t,y+h_j)$) is useful only when $(h_0,\ldots,h_d)$ are affinely independent
\begin{thm}\label{th:estuu}
Consider an admissible inviscid gas flow in $Q_T$, with finite total mass and energy. Then one has
\begin{equation}\label{eq:estuu}
\sup_{h_0,\ldots,h_d}\int_{Q_T}\left(\left(\prod_{j=0}^d\rho_j\right)\cdot({\rm Cor}(u_0,\ldots,u_d))^2\right)^{\frac1d}(t,y)\,dy\,dt\le_dM^{\frac1d}\sqrt{ME_0\,}\,.
\end{equation}
\end{thm}

We point out that the integrand in (\ref{eq:estuu}) shares the same physical dimension as the expected $\rho^{1+\frac1d}|u|^2$. Recall in passing that the product $ME_0$ in the right-hand side of (\ref{eq:estuu}) may be replaced by its Galilean invariant version, exactly like in \cite{Ser_DPT}~:
$$\frac14\int_{\R^d}\!\int_{\R^d}\rho(0,y)\rho(0,z)|u(0,z)-u(0,y)|^2dy\,dz+M\int_{\R^d}(\rho e)(0,y)\,dy,$$
where the last integral, the internal energy, remains unchanged.

Once again (\ref{eq:estuu}) has the flavor of a Strichartz inequality: denoting the space integral
$$H(t,h_0,\ldots,h_d):=\int_{\R^d}\left(\left(\prod_{j=0}^d\rho_j\right)\cdot({\rm Cor}(u_0,\ldots,u_d))^2\right)^{\frac1d}(t,y)\,dy,$$
it tells us that $H\in L^\infty_{h_0,\ldots,h_d}L^1_t$. This must be put in front of a direct application of (\ref{eq:totME}), which gives instead $H\in L^\infty_{t,h_1,\ldots,h_d}L^d_{h_0}\,$, see Paragraph \ref{ss:vel}. 

We observe, sadly, that our collection of {\em a priori} estimates, namely the conservation of mass, decay of energy, together with (\ref{eq:pGD}) and (\ref{eq:estuu}), is not strong enough to give a meaning to the flux in the balance law of energy (which is a conservation law in some cases):
\begin{equation}\label{eq:consen}
\partial_t\left(\frac12\,\rho|u|^2+\rho e\right)+{\rm div}_y\left[\left(\frac12\,\rho|u|^2+\rho e+p\right)u\right]\le0.
\end{equation}
As a matter of fact, we control quantities which either are quadratic in the velocity field, or do not involve it, while the flux in (\ref{eq:consen}) is cubic in $u$.

\bigskip

Another significant improvement of the corpus of estimates is the following statement concerning the pressure.
\begin{thm}\label{th:pconvol}
Admissible flows of Gas dynamics, with finite total mass and energy, satisfy
\begin{equation}
\label{eq:Schurp}
\sup_{\tau,\eta}\int_0^\infty\int_{\R^d}\left(\frac{(t-\tau)^2}{(E_0(t-\tau)^2+M|y-\eta|^2)^{\frac d2+1}}\right)^{\frac1d}p(t,y)dy\,dt\le_d E_0^{1-\frac1d}.
\end{equation}
\end{thm}
Inequality (\ref{eq:Schurp}) is reminiscent of a well-known fact about convolution of functions $f:\R^n\to\R$ of bounded variations, say with compact support,
\begin{equation}\label{eq:Luc}
\left\|\frac1{|x|}*f\right\|_{L^\infty}\le_nTV(f).
\end{equation}

\begin{rque}
\label{r:Luc}
According to L. Tartar (personal communication), the latter is a consequence of the improved Sobolev--Lorentz embedding (see A. Alvino \cite{Alv})
\begin{equation}
\label{eq:betterSob}
BV(\R^n)\subset L^{\frac n{n-1},1}(\R^n),
\end{equation} 
with the fact that $x\mapsto \frac1{|x|}$ belongs to the dual space $L^{n,\infty}(\R^n)$. Luc pointed out that the convolution product actually belongs to ${\cal C}_b(\R^n)$, the space of continuous bounded functions, because of the density of ${\cal C}_K(\R^n)$ in $L^{\frac n{n-1},1}(\R^n)$.
\end{rque}

 We shall give below a new proof of (\ref{eq:Luc}), by means of Compensated Integrability, and extend it to the convolution by arbitrary positively homogeneous kernels of degree $-1$, see Proposition \ref{p:convolg}.

\bigskip

\paragraph{Plan of the paper.} Section \ref{s:newCI} displays new forms of Compensated Integrability. The homogeneous polynomial $A\mapsto \det A$ being hyperbolic over ${\bf Sym}_n^+$ in the sense of G{\aa}rding \cite{Gaa}, we may use the properties of its multi-linearization, the {\em mixed determinant}. Part \ref{s:BV}, which applies the abstract results to scalar functions, is two-fold. It begins with BV-functions and continues with a Gagliardo-like inequality for time-dependent functions. Section \ref{s:newGD} treats applications to Gas dynamics, and contains the proofs of Theorems \ref{th:estuu} and \ref{th:pconvol}. In the same spirit, it displays a new estimate of the so-called {\em defect measure} introduced by \cite{LM} when averaging renormalized solutions of Boltzman equation.

\bigskip

\paragraph{Acknowledgement.} I am grateful to Luc Tartar for sharing his precious knowledge of the history of Sobolev--Lorentz embeddings, and to P. Mironescu for his interest in this work.

\section{New forms of Compensated Integrability}\label{s:newCI}

We start with a refined result, whose proof illustrates the scaling technique which is at stake almost everywhere in the theory.
\begin{thm}\label{th:prod}
Let $A\succ0_d$ be Div-BV over $\R^n$. Then 
\begin{equation}\label{eq:refCI}
\int_{\R^n}(\det A)^{\frac1{n-1}}dx\le n^{\frac n{n-1}}c_n\left(\prod_{i=1}^n\|({\rm Div}\,A)_i\|_{\cal M}\right)^{\frac1{n-1}}.
\end{equation}
\end{thm}

\bigskip

\bepr

We rescale both dependent and independent variables,
$$x_j'=\mu_jx_j,\qquad a_{ij}'=\mu_i\mu_ja_{ij}.$$
Since $A'=\Delta A\Delta$ for a diagonal matrix $\Delta$, it is still positive semi-definite.
Because $\partial_j'=\frac1{\mu_j}\partial_j$, we have
$$({\rm Div}'A')_i=\mu_i({\rm Div}\,A)_i$$
and thus $A'$ is also Div-BV. Let us apply Thm \ref{th:fund} to $A'$, using
$$\det A'=\left(\prod_{j=1}^n\mu_j\right)^2\det A,\qquad dx'=\left(\prod_{j=1}^n\mu_j\right)dx,$$
and
$$|{\rm Div}\,A|\le\sum_{j=1}^n|({\rm Div}\,A)_j|.$$
We obtain
$$\int_{\R^n}(\det A)^{\frac1{n-1}}dx\le c_n\left(\prod_{j=1}^n\mu_j\right)^{-\frac1{n-1}}\left(\sum_{i=1}^n\mu_i\|({\rm Div}\,A)_i\|_{\cal M}\right)^{\frac n{n-1}}.$$
Chosing
$$\mu_i=\frac1{\|({\rm Div}\,A)_i\|_{\cal M}}\,,$$
we obtain (\ref{eq:refCI}).

\enpr

\bigskip

\begin{rque} The inequality in Theorem \ref{th:prod} can be recast as
$$\log\int_{\R^n}(\det A)^{\frac1{n-1}}dx\le \log(n^{\frac n{n-1}}c_n)+\frac1{n-1}\,\sum_{j=1}^n\log\|({\rm Div}\,A)_j\|_{\cal M}.$$
Since it
remains valid when the canonical basis is replaced by an arbitrary orthonormal basis, we obtain by averaging
$$\log\int_{\R^n}(\det A)^{\frac1{n-1}}dx\le \log(n^{\frac n{n-1}}c_n)+\frac n{n-1}\,\Xint-_{S_{n-1}}\log\|({\rm Div}\,A)\cdot e\|_{\cal M}ds(e).$$
In other words, we have
\begin{equation}
\label{eq:multFI}
\int_{\R^n}(\det A)^{\frac1{n-1}}dx\le c_n\left(n\exp\Xint-_{S_{n-1}}\log\|{\rm div}(Ae)\|_{\cal M}ds(e)\right)^{\frac n{n-1}}.
\end{equation}
\end{rque}

\subsection{Estimating mixed determinant.}

To go further, we invoque the notion of {\em mixed determinant}. This is the symmetric $n$-linear form $D_n$ with the property that
$$D_n(M,\ldots,M)=\det M,\qquad\forall M\in{\bf M}_n(\R).$$
For instance, the case $n=2$ (where the determinant is a quadratic form) gives
$$D_2(M,M')=\frac12\,\left(\det(M+M')-\det M-\det M'\right)=\frac14\,\left(\det(M+M')-\det (M-M')\right).$$
More generally, we have
$$D_n(M_1,\ldots,M_n)=\frac1{2^n\,n!}\,\sum_{\epsilon_1,\ldots,\epsilon_n=\pm1}\epsilon_1\cdots\epsilon_n\det(\epsilon_1M_1+\cdots+\epsilon_nM_n).$$
We notice the useful instance
\begin{equation}
\label{eq:AAA}
D_n(B,M,\ldots,M)=\frac1n\,\Tr(B^T\widehat M).
\end{equation}

\paragraph{Hyperbolicity.}
Let us recall that the homogeneous polynomial $\det:{\bf Sym}_n\rightarrow\R$ is {\em hyperbolic} in G{\aa}rding's terminology \cite{Gaa}, meaning that there exists a direction $E\ne0_n$ such that $\det E>0$ and for every $A\in{\bf Sym}_n$, the roots of the univariate polynomial $t\mapsto\det(A-tE)$ are real.  The connected component $\Gamma$ of $E$ in the complement of $\{A|\det A=0\}$ is a convex cone, called the {\em forward cone}. Actually, the hyperbolic property stands in every direction $E'\in\Gamma$. We may of course choose $E=I_n$, and we see that $\Gamma={\bf SPD}_n$ is the positive definite cone.

G{\aa}rding proved that a hyperbolic polynomial satisfies the reverse H\"older inequality in the closure of  its forward cone. In the case of the determinant, this reads
\begin{prop}\label{p:Gaa}
For every $A_1,\ldots,A_n\in{\bf Sym}_n^+$, there holds
\begin{equation}
\label{eq:Gaa}
\left(\prod_{j=1}^n\det A_j\right)^{\frac1n}\le D_n(A_1,\ldots,A_n).
\end{equation}
\end{prop}
In particular the mixed determinant takes non-negative values over ${\bf Sym}_n^+\times\cdots\times{\bf Sym}_n^+$. Developing the multilinear expression $D_n(A_1+\cdots+A_n,\ldots,A_1+\cdots+A_n)$, we deduce
\begin{equation}
\label{eq:multi}
D_n(A_1,\ldots,A_n)\le\frac1{n!}\,\det(A_1+\cdots+A_n),\qquad\forall\,A_1,\ldots,A_n\in{\bf Sym}_n^+.
\end{equation}

Applying Theorem \ref{th:fund} to $A=A_1+\cdots+A_n$, and using (\ref{eq:multi}), we deduce an inequality for non-negative Div-BV tensors over $\R^n$~:
$$\int_{\R^n}(D_n(A_1,\ldots,A_n))^{\frac1{n-1}}dx\le_n\left(\sum_{j=1}^n\|{\rm Div}\,A_j\|_{\cal M}\right)^{\frac n{n-1}}.$$
Replacing $A_j$ by $\lambda_jA_j$ where $\lambda_j=\|{\rm Div}\,A_j\|_{\cal M}^{-1}$, we infer the estimate
\begin{thm}\label{th:mulest}
Given positive semi-definite Div-BV tensors $A_1,\ldots,A_n$ over $\R^n$, one has
\begin{equation}\label{eq:mulest}
\int_{\R^n}(D_n(A_1,\ldots,A_n))^{\frac1{n-1}}dx\le_n\left(\prod_{j=1}^n\|{\rm Div}\,A_j\|_{\cal M}\right)^{\frac 1{n-1}}.
\end{equation}
\end{thm}

\subsection{Estimating Schur complements}

Consider a positive semi-definite Div-BV tensor written blockwise
$$A=\begin{pmatrix} \rho & m^T \\ m & B \end{pmatrix},
$$
where $\rho$ is scalar. The positiveness tells us on the one hand that $\rho\ge0$ and on the other hand that either $\rho=0$, $m=0$ and $B\succ0_{n-1}$, or $\rho>0$ and $B\succ\frac{m\otimes m}\rho$\,. Thus let us assume that $\rho>0$, and denote $S:=B-\frac{m\otimes m}\rho\succ0_{n-1}$ its {\em Schur complement}. Then $A$ splits into the sum of two positive semi-definite matrices:
\begin{equation}
\label{eq:Aschur}
A=\begin{pmatrix} \rho & m^T \\ m & \frac{m\otimes m}\rho \end{pmatrix}+\begin{pmatrix} 0 & 0^T \\ 0 & S \end{pmatrix}=:\begin{pmatrix} \rho & m^T \\ m & \frac{m\otimes m}\rho \end{pmatrix}+A'.
\end{equation}
The {\em Schur complement formula}
\begin{equation}\label{eq:Scf}
\det A=\rho\det S
\end{equation}
allows us to make calculations without expressing $S$ in closed form.

\bigskip

Suppose now that $F\succ0_n$ is another Div-BV tensor. Developping the expression $D_n(F,A,\ldots,A)$ and discarding all non-negative terms but one, we obtain
\begin{equation}
\label{eq:FAA}
D_n(F,A,\ldots,A)\ge D_n(F,A',\ldots,A')=\frac1n\,\Tr(F\widehat{A'})=\frac1n\,f_{11}\det S.
\end{equation}
We infer the Schur complement estimate:
\begin{prop}\label{p:Sch}
Given two positive semi-definite Div-BV tensors $F,A$ over $\R^n$, with $a_{11}>0$,  the Schur complement $S$ of $a_{11}$ satisfies
\begin{equation}
\label{eq:Sch}
\int_{\R^n}(f_{11}\det S)^{\frac1{n-1}}dx\le_n\|{\rm Div}\,A\|_{\cal M}\|{\rm Div}\,F\|_{\cal M}^{\frac1{n-1}}.
\end{equation}
\end{prop}

\bigskip

\begin{rque}
\label{rk:SigSch}
When the entry $a_{11}$ is not strictly positive everywhere, it might sometimes be difficult, if not impossible, to define its Schur complement. Instead, it may happen that the tensor decomposes as $A=K+A'$ with $K\succ0_n$ and 
$$A'=\begin{pmatrix} 0 & 0 \\ 0 & \Sigma \end{pmatrix},\qquad \Sigma\succ0_{n-1}.$$
Then the same argument as above works: (\ref{eq:FAA}) still holds true, and implies a generalization of (\ref{eq:Sch}),
\begin{equation}
\label{eq:SigSch}
\int_{\R^n}(f_{11}\det \Sigma)^{\frac1{n-1}}dx\le_n\|{\rm Div}\,A\|_{\cal M}\|{\rm Div}\,F\|_{\cal M}^{\frac1{n-1}}.
\end{equation}
\end{rque}

\bigskip

The full strength of (\ref{eq:SigSch}) occurs when we choose an extreme  tensor $F$. By ``extreme'', we mean that $F$ be positively homogeneous of degree $1-n$, see \cite{Ser_JMPA}. Since the homogeneity is not compatible with $f_{ij}\in L^1(\R^n)$, we use a cut-off function $\phi(r)$ before getting rid of it:
$$F=\phi(r)G,\qquad G(x)=\frac{x\otimes x}{r^{n+1}}\,,\qquad r=|x|.$$
The cut-off is a non-increasing function such that $\phi(r)\equiv1$ over $(0,R)$ and $\equiv0$ over $(R+1,+\infty)$. Since ${\rm Div}\,G\equiv0$, we have 
$$\|{\rm Div}\,F\|_{\cal M}=\int_{\R^d}|G\nabla\phi|\,dx=\int_{\R^d}|\phi'(r)|\,\frac{dx}{r^{n-1}}=|S_{n-1}|\int_0^\infty|\phi'(r)|\,dr=|S_{n-1}|.$$
Applying (\ref{eq:Sch}), we obtain
$$\int_{|x|<R}\left(\frac{x_1^2}{r^{n+1}}\right)^{\frac1{n-1}}(\det \Sigma)^{\frac1{n-1}}dx\le_n\|{\rm Div}\,A\|_{\cal M}.$$
Taking the limit as $R\to+\infty$ and relaxing the position of the origin, we obtain the estimate
\begin{equation}
\label{eq:estunun}
\sup_{\xi\in\R^n}\int_{\R^n}\left(\frac{(x_1-\xi_1)^2}{|x-\xi|^{n+1}}\right)^{\frac1{n-1}}(\det \Sigma)^{\frac1{n-1}}dx\le_n\|{\rm Div}\,A\|_{\cal M}.
\end{equation}

When $A$ is positive definite, the Schur complement is well-defined. Then applying a rotation, using (\ref{eq:Scf}) and the Schur formula, we conclude
\begin{thm}\label{th:estSchur}
Let $A:\R^n\to{\bf SPD}_n$ be a Div-BV tensor. Then we have 
\begin{equation}
\label{eq:estSchur}
\sup_{\omega\in S_{n-1}}\sup_{\xi\in\R^n}\int_{\R^n}\left(\frac{(\omega\cdot(x-\xi))^2}{|x-\xi|^{n+1}}\right)^{\frac1{n-1}}\left(\frac{\det A}{\omega^TA\omega}\right)^{\frac1{n-1}}dx\le_n\|{\rm Div}\,A\|_{\cal M}.
\end{equation}
\end{thm}

\subsection{Estimating the rank-one part in (\protect\ref{eq:Aschur})}

Let us rewrite (\ref{eq:Aschur}) in the form $A=\rho U\otimes U+A'$, where
$$U=\binom1u,\qquad u=\frac m\rho\,.$$
Suppose that we are given a collection $(A_1,\ldots,A_n)$ of positive semi-definite Div-BV tensors over $\R^n$. 
Let us form the, still positive and Div-BV, tensor
$$\A=A_1+\cdots+A_n.$$
Decomposing $A_j=\rho_jU_j\otimes U_j+A_j'$, we have
$$\A\succ\rho_1U_1\otimes U_1+\cdots+\rho_nU_n\otimes U_n,$$
whence
$$\det\A\ge\det(\rho_1U_1\otimes U_1+\cdots+\rho_nU_n\otimes U_n).$$
Expressing
$$\rho_1U_1\otimes U_1+\cdots+\rho_nU_n\otimes U_n=U^TRU,\qquad R={\rm diag}(\rho_1,\ldots,\rho_n),$$
where $U$ stands for the matrix whose columns are $U_1,\ldots,U_n$, we infer
$$\det\A\ge(\det R)\cdot(\det U)^2=\left(\prod_{j=1}^n\rho_j \right)\cdot\left[{\rm Cor}(u_1,\ldots,u_n)\right]^2.$$

Assembling the material above, and applying a standard scaling argument, we end up with
\begin{prop}\label{p:AunAn}
Let $A_1,\ldots,A_n$ be positive semi-definite Div-BV tensors over $\R^n$, written blockwise 
$$A_j=\rho_j\binom1{u_j}\otimes\binom1{u_j}+\begin{pmatrix} 0 & 0 \\ 0 & B_j \end{pmatrix}.$$ 
Then
$$\int_{\R^n}\left(\prod_{j=1}^n\rho_j\cdot\left[{\rm Cor}(u_1,\ldots,u_n)\right]^2 \right)^{\frac1{n-1}}dx\le_n\left(\sum_{j=1}^n\|{\rm Div}\,A_j\|_{\cal M}\right)^{\frac n{n-1}}.$$
\end{prop}
Remark that the integrand in the left-hand side has the dimension of a $\rho^{\frac n{n-1}}|u|^2$.

\section{Applications to scalar functions}\label{s:BV}

\subsection{BV functions}

If $f\in BV(\R^n)$ is non-negative, we may apply (\ref{eq:mulest}) to the choices
$$A_1=F,\qquad A_2=\cdots=A_n=fI_n.$$
With (\ref{eq:AAA}), we find
$$D_n(F,fI_n,\ldots,fI_n)=\frac1n\,f^{n-1}\Tr F=f^{n-1}\frac{\phi}{nr^{n-1}}\,.$$
With ${\rm Div}\,(fI_n)=\nabla f$, we thus obtain
$$\int_{\R^n}f(x)\phi(r)^{\frac1{n-1}}\,\frac{dx}r\le_n\|{\rm Div}\,(fI_n)\|_{\cal M}= TV(f).$$
Letting $R\to+\infty$, this yields
$$\int_{\R^n}f(x)\,\frac{dx}r\le_n TV(f).$$
Eventually, placing the origin at an arbitrary point, this rewrites
$$\forall\xi\in\R^n,\qquad\int_{\R^n}f(x)\,\frac{dx}{|x-\xi|}\le_n TV(f).$$
To pass from non-negative functions to signed functions, we may use $|r^{-1}*f|\le r^{-1}*|f|$, and $TV(|f|)\le TV(f)$. This re-proves (\ref{eq:Luc}).

\bigskip

The calculation above can be generalized as follows. If $g\in L^{n-1}(S_{n-1})$ is non-negative, then choose instead
$$\bar G(x)=g\left(\frac xr\right)^{n-1}\frac{x\otimes x}{r^{n+1}}$$
and $F(x)=\phi(r)\bar G(x)$ as before. In \cite{Ser_JMPA} we have seen that ${\rm Div}\,\bar G=V\delta_{x=0}$ where 
$$V=\int_{S_{n-1}}g(\omega)^{n-1}\omega\,ds(\omega).$$
With $\phi$ as before, this yields
$${\rm Div}\,F=\phi(0)V\delta_{x=0}+g\left(\frac xr\right)^{n-1}\frac{\phi'}{r^{n-1}}\,\vec e_r,$$
Whence
$$\|{\rm Div}\,F\|_{\cal M}\le2\|g\|_{L^{n-1}(S_{n-1})}^{n-1}.$$
Then (\ref{eq:mulest}), applied to $(F,fI_n,\ldots,fI_n)$ yields
$$\int_{\R^n}f(x)\phi(r)^{\frac1{n-1}}g\left(\frac xr\right)\,\frac{dx}r\le_n TV(f)\|g\|_{L^{n-1}(S_{n-1})}.$$
Taking as above the supremum over the admissible $\phi$ (the limit as $R\to+\infty$), and chosing arbitrarily the origin, we end up with (see also Remark \ref{r:Luc} for the replacement of $L^\infty$ by ${\cal C}_b$)
\begin{prop}\label{p:convolg}
There exists a finite constant $C_n$ such that, for every $g\in L^{n-1}(S_{n-1})$, there holds 
$$\|\bar g*f\|_{L^\infty}\le C_n\cdot TV(f)\|g\|_{L^{n-1}(S_{n-1})},\qquad \forall \,f\in BV(\R^n),\qquad\bar g(x):=\frac1r\,g\left(\frac xr\right)\,.$$
\end{prop}

\subsection{A Gagliardo-like inequality}

Let us recall Gagliardo's inequality \cite{Gag} in $\R^d$~: given $d$ functions $f_j\in L^{d-1}(\R^{d-1})$, the new function 
$$f(y)=\prod_{j=1}^df_j(\hat y_j),\qquad y\in\R^d$$
is integrable. And there is a functional inequality
$$\|f\|_{L^1(\R^d)}\le\prod_{j=1}^d\|f_j\|_{L^{d-1}(\R^{d-1})}.$$
If $d=2$, this is nothing but Fubini's theorem. For general $d$, this can be recovered, up to a multiplicative constant, by applying (\ref{eq:refCI}) to the tensor $\phi\,{\rm diag}(|f_1|^{d-1},\ldots,|f_d|^{d-1})$, where $\phi$ is a cut-off as above, see \cite{Ser_DPT}.

\bigskip

We now introduce a time parameter, thus having $d$ functions $f_j(t,\hat y_j)$ and defining as above
\begin{equation}\label{eq:alaGag}
f(t,y)=\prod_{j=1}^df_j(t,\hat y_j),\qquad (t,y)\in\R^{1+d}.
\end{equation}
Mind that each $f_j$ is defined over a $d$-dimensional space. Let us set $n=1+d$ and form the tensor
$$A(x)=\frac{\phi(r) x\otimes x}{r^{n+1}}+{\rm diag}\left(0,\phi(|y_1|)|f_1(\hat x_1)|^{n-1},\ldots,\phi(|y_d|)f_d(\hat x_d)^{n-1}\right)=:\sum_{i=0}^dA_i.$$
We have as above
$$\|{\rm Div}\,A_0\|_{\cal M}=|S_{n-1}|,$$
while 
$$\|{\rm Div}\,A_j\|_{\cal M}=\int_{\R^n}|\phi'(|y_j|)|\cdot|f_j(\hat x_j)|^{n-1}dx=2\|f_j\|_{L^{d}(\R^{d})}^d.$$
By Schur formula, 
$$\det A=\frac{\phi(r)t^2}{r^{n+1}}\cdot\prod_{j=1}^d\phi(|y_j|)|f_1(\hat x_j)|^{n-1}.$$
Applying Theorem \ref{th:fund}, we infer
$$\int_{|x|<R}\left(\frac{t^2}{(t^2+|y|^2)^{1+\frac d2}}\right)^{\frac1d}\left|\prod_{j=1}^df_j(t,\hat y_j)\right|\,dy\,dt\le_d\left(1+\sum_{j=1}^d\|f_j\|_{L^d(\R\times\R^{d-1})}^d\right)^{1+\frac1d}.$$
Passing to the limit as $R\to+\infty$, and making the appropriate scaling, we conclude
\begin{thm}\label{th:alaGag}
Given the functions $f_1,\ldots,f_d\in L^d(\R\times\R^{d-1})$, the function $f$ defined over $\R\times\R^d$ by (\ref{eq:alaGag}) satisfies
$$\int_{\R\times\R^d}\left(\frac{t^2}{(t^2+|y|^2)^{1+\frac d2}}\right)^{\frac1d}|f(t,y)|\,dy\,dt\le_d\prod_{j=1}^d\|f_j\|_{L^{d}(\R\times\R^{d-1})}.$$
\end{thm}

\section{Applications to Gas dynamics}\label{s:newGD}

Once again we set $n=1+d$, where $d$ is the space dimension, so that $\frac1{n-1}=\frac1d$\,. The Euler system governs the conservation of mass and linear momentum
\begin{eqnarray*}
\partial_t\rho+{\rm div}_y(\rho u) & = & 0, \\
\partial_t(\rho u)+{\rm Div}_y(\rho u\otimes u)+\nabla_yp & = & 0. 
\end{eqnarray*}
This is recast as ${\rm Div}_{t,y}A=0$, where 
$$A=\begin{pmatrix} \rho & \rho u^T \\ \rho u & \rho u\otimes u+pI_d \end{pmatrix}=\rho U\otimes U+pJ\qquad\hbox{in }Q_T,$$
with 
$$U:=\begin{pmatrix} 1 \\ u \end{pmatrix},\qquad J=\begin{pmatrix} 0 & 0 \\ 0 & I_d \end{pmatrix}.$$
We recognize the Schur complement $pI_d$ of $\rho$ in $A$, whose determinant is just $p^d$. 

\subsection{A new estimate of the pressure}

Applying (\ref{eq:estunun}), we have
\begin{equation}\label{eq:nonhom}
\sup_{\tau,\eta}\int_{Q_T}\left(\frac{(t-\tau)^2}{((t-\tau)^2+|y-\eta|^2)^{\frac d2+1}}\right)^{\frac1d}p(t,y)dy\,dt\le_d M+\sqrt{ME_0}.
\end{equation}
Since the bound does not depend upon the length of the time interval, this estimate is valid over $Q_\infty$. 

The estimate above is not homogeneous from a Physical point of view. To prove Theorem \ref{th:pconvol}, we apply a scaling technique as in \cite{Ser_DPT}~: for every parameter $\mu>0$, the dependent/independent variables
$$t'=t,\quad y'=\mu y,\quad\rho'=\rho,\quad u'=\mu u,\quad p'=\mu^2 p$$
define an admissible flow, to which we may apply (\ref{eq:nonhom}). Chosing $\mu=\sqrt{M/E_0\,}$\,, we obtain (\ref{eq:Schurp}).

\bigskip

Notice that we may rewrite (\ref{eq:Schurp}) in terms of the square root mean velocity $\bar u:=\sqrt{2E_0/M\,}\,$, 
$$\sup_{\tau,\eta}\int_{Q_T}\left(\frac{(t-\tau)^2}{\left((t-\tau)^2+\frac{|y-\eta|^2}{\bar u^2}\right)^{\frac d2+1}}\right)^{\frac1d}p(t,y)dy\,dt\le_d E_0.$$

%

\subsection{The velocity field}\label{ss:vel}

To get more information about the velocity field than just what is given by the energy estimate, we apply Proposition \ref{p:AunAn} to the tensors $A_0,\ldots,A_d$ obtained from $A$ by a space shift:
$$A_j(t,y):=A(t,y+h_j).$$
Hereabove $h_0,\ldots,h_d$ are constant vectors. Of course the $A_j$'s are Div-BV whenever $A$ is so. We thus obtain the inequality
$$\int_{Q_T}\left(\prod_{j=0}^d\rho(t,y+h_j)\cdot\left({\rm Cor}(u(t,y+h_0),\ldots,u(t,y+h_d))\right)^2\right)^{\frac1d}dy\,dt\le_d\left(M+\sqrt{ME_0\,}\right)^{1+\frac1d}.$$
After the same scaling procedure as above, this gives Theorem \ref{th:estuu}.

\paragraph{Discussion.} The integral depends only upon $(h_0,\ldots,h_d)$ modulo translations~; we may thus assume $h_0=0$ without loss of generality.
Estimate (\ref{eq:estuu}) is a statement about the expression
$$H(t,h_1,\ldots,h_d)=\int_{\R^d}\left(\rho\prod_{j=1}^d\rho_j\cdot\left({\rm Cor}(u,u_1,\ldots,u_d)\right)^2\right)^{\frac1d}dy.$$
Namely, it says that
\begin{equation}
\label{eq:ptells}
H\in L^\infty_{h_1,\ldots,h_d}L^1_t.
\end{equation}
We shall compare (\ref{eq:ptells}) with what can be said of $H$ by using only the conservation of total mass and the decay of total energy. On the one hand, (\ref{eq:estuu}) is established through Compensated Integrability and thus requires that the mass-momentum of the fluid be positive semi-definite. As such, it applies to an inviscid fluid, as well as to the conservation of mass/momentum satisfied by the renormalized solutions of the Boltzman equation (see \cite{LM}), but it does not apply to a viscous fluid. On the other hand, the pure mass-energy estimate below is valid in a much more general context, since it does not assume the positiveness of the mass-momentum tensor.

For this direct estimate, we first notice that, $U=(U_0|\cdots|U_d)$ being the same matrix as above,
$$|\det U|=\left|\sum_{k=0}^d(-1)^k\det(\ldots,u_{k-1},u_{k+1},\ldots)\right|\le\sum_{k=0}^d\prod_{j\ne k}|u_j|.$$
We infer a majorization $H\le_d\sum_{k=0}^d H_k$ where
$$H_k(t,h_1,\ldots,h_d)=\int_{\R^d}\left(\prod_{j=0}^d\rho_j\cdot\prod_{j\ne k}|u_j|^2\right)^{\frac1d}dy.$$
The treatment of $H_k$ depends on whether $k=d$ or $k<d$, although it yields the same bound. Let us begin with the latter case. Applying the H\"older inequality with exponents $(d,\ldots,d)$, we have 
\begin{eqnarray*}
H_k^d & \le & \left(\prod_{j\ne k,d}\int_{\R^d}\rho_j|u_j|^2dy\right)\cdot\int_{\R^d}\rho_k(t,y)(\rho|u|^2)(t,y+h_d)\,dy \\
& \le & (2E_0)^{d-1}\int_{\R^d}\rho_k(t,y)(\rho|u|^2)(t,y+h_d)\,dy.
\end{eqnarray*}
Integrating in $h_d$ and using Fubini, this gives
$$\int_{\R^d}H_k(t,h_1,\ldots,h_d)^ddh_d\le_d (2E_0)^{d}\int_{\R^d}\rho_k(t,y)\,dy=M(2E_0)^d.$$
If instead $k=d$, we have
\begin{eqnarray*}
H_d^d & \le & \left(\prod_{j=1}^{d-1}\int_{\R^d}\rho_j|u_j|^2dy\right)\cdot\int_{\R^d}\rho(t,y+h_d)(\rho|u|^2)(t,y)\,dy \\
& \le & (2E_0)^{d-1}\int_{\R^d}\rho(t,y+h_d)(\rho|u|^2)(t,y)\,dy,
\end{eqnarray*}
and again
$$\int_{\R^d}H_d(t,h_1,\ldots,h_d)^ddh_d\le_d ME_0^d.$$
The {\bf direct} mass-energy bound for $H$ is thus
\begin{equation}
\label{eq:estdirect}
\sup_{t,h_1,\ldots,h_{d-1}}\int_{\R^d}H(t,h_1,\ldots,h_d)^ddh_d\le_d ME_0^d.
\end{equation}
Thus the mass-energy conservation (or decay) provides the qualitative result that
\begin{equation}
\label{eq:direct}
H\in L^\infty_{t,h_1,\ldots,h_{d-1}}L^d_{h_d}.
\end{equation}
Comparing with (\ref{eq:ptells}), this amounts to replacing $L^\infty_{h_d}L^1_t$ by $L^\infty_tL^d_{h_d}$. 

Somehow the inequality (\ref{eq:estdirect}) is brutal, in the sense that it applies to every quantity ${\cal H}$ of the form
$${\cal H}(t,h_1,\ldots,h_d)=\int_{\R^d}\left(\prod_{j=0}^d\rho_j\cdot F(u_0,\ldots,u_d)\right)^{\frac1d}dy$$
where $F$ is a homogeneous polynomial of global degree $2d$, quadratic in each argument.

\subsection{Renormalized solutions of Boltzman equation}

Let us consider  the Boltzman equation
\begin{equation}\label{kineq}
(\partial_t+v\cdot\nabla_y)f(t,y,v)=Q[f(t,y,\cdot)], \qquad (t,y)\in Q_T=(0,T)\times\R^d,\, v\in\R^d.
\end{equation}
The Cauchy problem consists in finding a solution which fits an initial data $f(0,y,v)=f_0(y,v)$. 
The existence of distributional solutions to this problem is not known, except in space dimension $d=1$. When $d\ge1$, R. DiPerna \& P.-L. Lions \cite{DPL} proved instead the existence of a weaker notion of solutions, called {\em renormalized}. We shall not give a precise definition of this notion, and we content ourselves to recall that it implies, at the macroscopic level, the conservation of mass and a weak form of the conservation of momentum, in the sense that 
\begin{eqnarray}
\label{Euldef}
\partial_t\rho+{\rm div}_ym=0, & & \partial_tm+{\rm Div}_y\left(\int_{\R^d}fv\otimes v\,dv+\Sigma\right)=0 \\
\nonumber
\rho(t,y):=\int_{\R^d}f\,dv, & & m(t,y):=\int_{\R^d}fv\,dv,\qquad \Sigma\succ0_d.
\end{eqnarray}
The quantities $\rho,m$ are the mass density and linear momentum. Compared to what is formally expected, the second equation above contains an additional term $\Sigma$, called the {\em defect measure}, which takes values in ${\bf Sym}_d^+$~; see \cite{LM}. Finally, it is known that the total mass
$$M=\int_{\R^d}\int_{\R^d}f_0(y,v)\,dv\,dy\equiv\int_{\R^d}\int_{\R^d}f(t,y,v)\,dv\,dy$$
is a constant of the motion, and the total energy
$$E(t):=\frac12\int_{\R^d}\!\!\int_{\R^d}f(t,y,v)|v|^2dv\,dy+\frac12\,\int d\Tr\Sigma(t,\cdot)$$
is a non-increasing function of time and satisfies
$$E(t)\le E_0:=\frac12\int_{\R^d}\!\!\int_{\R^d}f_0(y,v)|v|^2dv\,dy.$$

The equations (\ref{Euldef}) can be recast by saying that the following tensor 
$$A=\begin{pmatrix} \int_{\R^d}f\,dv & \int_{\R^d}fv^T\,dv \\ \int_{\R^d}fv\,dv & \int_{\R^d}fv\otimes v\,dv+\Sigma \end{pmatrix}$$
is Div-free.
For finite total mass and energy, the extension of $A$ by $0_n$ away from $Q_T=(0,T)\times\R^d$ is Div-BV. 

As mentionned in Remark \ref{rk:SigSch}, it might be difficult to construct the Schur complement of the density, if it vanishes here and there. Instead, decomposing $A=B+A'$ where
$$B(t,y)=\int_{\R^d}f(t,y,v)\binom1v\otimes\binom1v\,dv,\qquad A'=\begin{pmatrix} 0 & 0 \\ 0 & \Sigma \end{pmatrix},$$
we may apply (\ref{eq:SigSch}) to obtain
$$\sup_{\tau,\eta}\int_0^T\int_{\R^d}\left(\frac{(t-\tau)^2}{((t-\tau)^2+|y-\eta|^2)^{1+d/2}}\right)^{\frac1d}d(\det\Sigma)^{\frac1d}\le_dM+\sqrt{ME_0\,}\,,$$
where we recall that $(\det\Sigma)^{\frac1d}$ is a finite measure, satisfying
$$0\le(\det\Sigma)^{\frac1d}\le\frac1d\,\Tr\Sigma.$$
Letting $T\to+\infty$ (as usual, the upper bound does not depend upon the length of the time interval), and applying the usual scaling trick, we conclude that the defect measure is constrained by
\begin{equation}
\label{eq:DFmeas}
\sup_{\tau,\eta}\int_0^T\int_{\R^d}\left(\frac{(t-\tau)^2}{(E_0(t-\tau)^2+M|y-\eta|^2)^{1+d/2}}\right)^{\frac1d}d(\det\Sigma)^{\frac1d}\le_dE_0^{1-\frac1d}.
\end{equation}
Inequality (\ref{eq:DFmeas}) tells us that $(\det\Sigma)^{\frac1d}$ is not too singular. For instance it does not charge points (no Dirac mass)~; $\Sigma$ itself might charge points, but its density at a Dirac mass must be a singular matrix.

\begin{english}

\end{english}

\end{document}